\documentclass[preprint,12pt]{elsarticle}
\usepackage[top=2.5cm, bottom=2.5cm, left=2.5cm, right=2cm]{geometry}
\usepackage{amssymb,latexsym,amscd,amsmath,amsfonts,enumerate,supertabular}
\usepackage{graphicx}
\usepackage{tabularx}
\usepackage{subfigure}
\usepackage{enumerate}
 \usepackage{epstopdf}
 \journal{Elsevier}
\newtheorem{theorem}{Theorem}
\newtheorem{lemma}{Lemma}
\newtheorem{proposition}{Proposition}

\newtheorem{remark}{Remark}

\newenvironment{proof}[1][Proof]{\begin{trivlist}
\item[\hskip \labelsep {\bfseries #1}]}{\end{trivlist}}

\begin{document}
\begin{frontmatter}
\title{{\bf Nonstandard finite difference schemes\\ for a general predator-prey system}}

%\title{{\bf Transform a general continuous-time predator–prey model into a dynamically consistent discrete-time model by using nonstandard finite difference method}}
\author[1]{Quang A Dang}
\ead{dangquanga@cic.vast.vn}
\author[2]{Manh Tuan Hoang}
\ead{hmtuan01121990@gmail.com}
\address[1]{Center for Informatics and Computing, Vietnam Academy of Science and Technology (VAST),\\ 
18 Hoang Quoc Viet, Cau Giay, Hanoi, Vietnam}
\address[2]{Institute of Information Technology, 
Vietnam Academy of Science and Technology (VAST), \\
18 Hoang Quoc Viet, Cau Giay, Hanoi, Vietnam}
 \begin{abstract}
\small
In this paper we transform a continuous-time predator-prey system with general functional response and recruitment for both species into a discrete-time model by nonstandard finite difference scheme (NSFD). The NSFD model shows complete dynamic consistency with its continuous counterpart for any step size. 
Especially, the global stability of a non-hyperbolic equilibrium point in a particular case of parameters is proved by the Lyapunov stability theorem. The performed numerical simulations confirmed the validity of the obtained theoretical results.
\small
 \end{abstract}
\small 
\small  
\end{frontmatter}
\section {Introduction}
Predator-prey systems are among the most interested topics in mathematical biology and ecology. Their dynamics continues to attract attention from both applied mathematicians and ecologists because of its universal existence and importance \cite{Linda, Berryman, Brauer, Keshet}. The majority of the results in this topic is mainly concentrated on the study of the qualitative aspects of the continuous systems described by  systems of differential equations. It is worth to mention some recent typical works such as \cite{He, Hu, Lan, Qi, Shi, Wang}\ldots In that time the conversion of the continuous systems to discrete systems preserving the properties of the original continuous systems is of great importance. This problem has attracted attention from many researchers and a lot of results are reached for some dynamical systems. Nevertheless, some continuous models are studied fully from theoretical point of view but their discrete counterparts still are not investigated.\par

Biological systems including predator-prey systems often are described by ordinary or partial differential equations. There are many ways for converting continuous models to discrete counterparts. The most popular way for this purpose is to use standard difference methods such as Euler, Runge-Kutta schemes. However, in many nonlinear problems the standard difference schemes reveal a serious drawback which is called "numerical instabilities" \cite{Mickens1, Mickens2, Mickens4}. Under this concept we have in mind the phenomena when the discrete models, for example, the difference schemes, do not preserve properties of the corresponding differential equations. In \cite{Mickens1, Mickens2,Mickens3, Mickens4} Mickens showed many examples and analysed the numerical instabilities when using standard difference schemes. In general, standard difference schemes preserve the properties of the differential equations only in the case if the discretization parameter $h$ is sufficiently small. Therefore, when studying dynamical models in large time intervals the selection of small time steps will requires very large computational effort, so these discrete models are inefficient. Besides, for some special dynamical problems standard difference schemes cannot preserve the properties of the problems for any step sizes.\par

In order to overcome the numerical instabilities phenomena in 1989  Mickens \cite{Mickens0} introduced the concept \emph{Nonstandard Finite Difference} (NSDF) schemes and after that has developed NSDF methods in many works, such as \cite{Mickens1, Mickens2,Mickens3, Mickens4}. According to Mickens, NSDF schemes are those constructed following a set of five basic rules. The NSDF schemes preserve main properties  of the differential counterparts, such as positivity, monotonicity, periodicity, stability and some other invariants including energy and geometrical shapes. It should be emphasized that NSFD schemes can preserve all properties of the continuous models for any discretization parameters. The discrete models with these properties are called \emph{dynamically consistent}  \cite{AL1, AL2, DQA, DK1, DK2,DK3,DK4,DK5, Mickens4, Mickens5, Roeger3}.\par

Up to date NSFD schemes became a power and efficient tool for simulating dynamical systems, especially in converting continuous models to dynamically consistent discrete counterparts \cite{AL1, AL2, Cresson, Darti, DK1, DK2,DK3,DK4, Mickens3, Partidar, Roeger4, Roeger5, Wood}. The majority of these models are met in physics, mechanics, chemistry, biology, epidemiology, finance, \ldots with complicated dynamical behaviour.
For predator-prey systems, some NSFD schemes, which are dynamically consistent with them, are constructed, such as \cite{Bairagi, Darti, DK3, DK5}. Below we mention some models.

%%%%%%%%%%%%%%%%%%%%%%%%%%%%%%%%%%%%%%%%%%%%%%%%%%%%%%%%%%%%%%%%%%%%%%%%%%%%%%%
In 2006 Dimitrov and Kojouharov \cite{DK3} used NSFD schemes for constructing  discrete model for the general Rosenzweig-MacArthur predator-prey model with a logistic intrinsic growth of the prey population. The model has the form
\begin{equation*}
\dfrac{dx}{dt} = bx(x - 1) - a g(x)xy, \qquad \dfrac{dy}{dt} = g(x)xy - dy,
\end{equation*}
where $x$ and $y$ represent the prey and predator population sizes, respectively, $b > 0$ represents the intrinsic growth rate of the prey, $a > 0$ stands for the capturing rate and $d > 0$ is the predator death rate. Under some assumptions for the function  $g(x)$  the obtained discrete model preserves the positivity of solutions and locally asymptotical stability of the set of equilibrium points. The finite difference method is called positive and elementary stable nonstandard (PESN) method.\par
Next, in 2008  Dimitrov and Kojouharov \cite{DK5} constructed a PESN scheme for predator-prey models with general functional response of the form
\begin{equation*}
\dfrac{dx}{dt} = p(x) - af(x, y)y,  \qquad \dfrac{dy}{dt} = f(x, y)y - \mu(y),
\end{equation*}
where functions $p(x)$ and $\mu(y)$ describe the intrinsic growth rate of the prey and the mortality rate of the predator, respectively, the function $f(x, y)$ is called ''functional response'' and represents the per capita predator ''feeding rate'' per unit time.  Recently, in 2016, Bairagi and Biswas \cite{Bairagi} used PESN scheme for converting a predator-prey model with Beddington-DeAngelis Functional Response  to a dynamically consistent discrete model. The model under consideration is of the form
\begin{equation*}
\dfrac{dx}{dt} = x(x - 1) - \dfrac{\alpha x y}{1 + \beta x + \mu y}, \qquad \dfrac{dy}{dt} = \dfrac{E x y}{1 + \beta x + \mu y} - Dy,
\end{equation*}
where $\alpha, \beta, \mu, E, D$ are positive constants. This model is a particular case of the model considered by  Dimitrov and Kojouharov in \cite{DK5}.
 An another interesting predator-prey model, which should be mentioned is the harvesting Leslie-Gower predator-prey model. In the recent paper \cite{Darti} the authors also constructed NSFD scheme preserving the positivity of solutions and asymptotical stability for this model.\par
It should be emphasized that all the equilibrium points in the mentioned above predator-prey models are hyperbolic, therefore, for establishing the stability, it suffices to consider eigenvalues of the Jacobian matrix of the linearized system around equilibrium points.\par
Aiming at the conversion of continuous  predator-prey  models studied fully from the theoretical point of view to dynamically consistent discrete models, in the present paper  we consider a mathematical model for a predator-prey system with general functional response and recruitment, which also includes capture on both species \cite{Lindano}. The results of qualitative aspects of the model are given in  \cite{Lindano}. One important difference of this model in comparison with the above models is that the model \cite{Lindano} has one non-hyperbolic equilibrium point.\par

The paper is organized as follows. In Section 2 we recall from \cite{Lindano} the predator-prey system under consideration with the theoretical results of the existence of equilibrium points and their stability properties. In Section 3 we propose  NSFD schemes and study their positivity and existence of equilibrium points. Next, Section 4 is devoted to the stability analysis of the equilibrium points. Some numerical simulations for demonstrating the validity of the theoretical results obtained in the previous section are given in Section 5. Finally, Section 6 is Conclusion.
%%%%%%%%%%%%%%%%%%%%%%%%%%%%%%%%%%%%%%%%%%%%%%%%%%%%%%%%%%%%%%%%%%%%%%%
\section{Mathematical model}
We consider a mathematical model for a predator-prey system with general functional response and recruitment for both species \cite{Lindano}. The model is described by the system of nonlinear differential equations
\begin{equation}\label{eq:1}
\begin{split}
\dot{x}(t) &= x(t)f(x(t), y(t)) = x(t)\big[r(x(t)) - y(t)\phi(x(t)) - m_1 \big],\\
\dot{y}(t)& = y(t)g(x(t), y(t)) = y(t)\big[s(y(t)) + cx(t)\phi(x(t)) - m_2\big],
\end{split}
\end{equation}
where:
\begin{itemize}
\item   $x(t)$ and $y(t)$ are prey population and predator population, respectively; 
\item   $r(x)$ and $s(y)$ are the per capita recruitment rates of prey and predators, respectively; 
\item  $xy\phi (x)$ is the predator response, and $x\phi (x)$ is the number of prey consumed per predator in a unit of time;
\item  $c$ is a constant named conversion efficiency of prey into predators, generally $0 < c < 1$ , and $cxy\phi (x)$ is the predator numerical
response;
\item   $m_1 >0$ and $m_2>$ are the total mortality rates of prey and predators, respectively.
\end{itemize}

From the biological significance we have
\begin{equation}\label{eq:2}
\begin{split}
&\forall x \geq 0, \qquad r(x) > 0, \qquad r'(x) < 0, \qquad [xr(x)]' \geq 0, \qquad \text{ and } \qquad \lim_{x \to \infty}r(x) = 0,\\
&\forall y \geq 0, \qquad s(y) > 0, \qquad s'(y) < 0, \qquad [ys(y)]' \geq 0, \qquad \text{ and } \qquad \lim_{y \to \infty}s(y) = 0,\\
&\forall x \geq 0, \qquad \phi(x) > 0, \qquad \phi'(x) \leq 0, \qquad \text{and } \qquad [x\phi(x)]' \geq 0.
\end{split}
\end{equation}
It is easy to check that the region $\Omega = \mathbb{R}_+^2$ is a positive invariant set for the system \eqref{eq:1}. Other qualitative properties of the model \eqref{eq:1} including the existence of equilibrium points, their stability, are studied in detail in \cite{Lindano}. \\ 
In order to easily track the paper, below we mention these properties.
\begin{theorem}\label{Propositionp1m}(Existence of equilibrium points) \cite[Proposition 1]{Lindano}\\
System \eqref{eq:1} has four distinct kinds of possible equilibrium points in $\Omega$:
\begin{enumerate}[(i)]%for small alpha-characters within brackets.
\item A trivial equilibrium point $P_0^* = (x_0^*, y_0^*) = (0, 0)$, for all the values of the parameter.
\item An equilibrium point of the form $P_1^* = (x_1^*, y_1^*) = (K, 0)$, with $r(K) = m_1$, if and only if $m_1 < r(0)$.
\item An equilibrium point of the form $P_2^* = (x_2^*, y_2^*) = (0, M)$, with $s(M) = m_2$, if and only if $m_2 < s(0)$.
\item An equilibrium point of the form $P_3^* = (x_3^*, y_3^*) = (x^*, y^*)$, where $x^*$ satisfies the equation
\begin{equation*}
cx^*\phi(x^*) + s\Big(\dfrac{r(x^*) - m_1}{\phi(x^*)}\Big) - m_2 = 0,
\end{equation*}
and $y^*$ is given, as a function of $x^*$, by
\begin{equation*}
y^* = \dfrac{r(x^*) - m_1}{\phi(x^*)},
\end{equation*}
if and only if $(m_1, m_2)$ verifies $m_1 < r(0) - M\phi(0)$ and $m_2 < s(0)$ or $m_1 < r(0)$ and $s(0) < m_2 < s(0) + cK\phi(K)$.
\end{enumerate}
\end{theorem}
%%%%%%%%%%%%%%%%%%%%%%%%%%%%%%%%%%%%%%%%%%%%%%%%%%%%%%%%%%%%%%%%%%%%%%%%%%%%%%
\begin{theorem}\label{Propositionp2m}(Stability analysis)
\begin{enumerate}[(i)]%for small alpha-characters within brackets.
\item If $m_1 > r(0)$ and $m_2 > s(0)$, then the extinction equilibrium point $P_0^* = (0, 0)$ is locally asymptotically stable, and unstable otherwise.
\item If $m_1 \geq r(0)$ and $m_2 \geq s(0)$ then the extinction equilibrium point $P_0^*$ is globally asymptotically stable.
\item If $m_1 < r(0)$ and $m_2 > s(0) + cK\phi(K)$, then the equilibrium point of the form $P_1^* = (K, 0)$ is locally asymptotically stable, and unstable otherwise.The equilibrium point $(K, 0)$ shall be called the equilibrium point of extinction of the predator species.
\item If $m_1 > r(0) - M\phi(0)$ and $m_2 < s(0)$, then the equilibrium point of the form $P_2^* = (0, M)$ is locally asymptotically stable, and unstable otherwise.
The equilibrium point $(0, M)$ shall be called the equilibrium point of extinction of the prey species.
\item If an equilibrium point of the form $P_3^* = (x^*, y^*)$ belongs to $\Omega$, then it is locally asymptotically stable. This equilibrium point $(x^*, y^*)$ shall be named the ecological stability equilibrium.
\end{enumerate}
\end{theorem}
In general, the property of stability of the set of equilibria of differential equations  plays an essential role in the study of asymptotical behavior of the solutions of differential equations. The construction of difference schemes, which preserve the stability of the equilibrium points, is important in numerical simulation of differential equations. The difference schemes with this stability property is called \emph{elementary stable} schemes  \cite{AL1, DK1, DK2, Wood}. There are many works concerning the elementary stable schemes. The typical results are for general dynamical systems \cite{DK1, DK2, K1} and for other specific systems \cite{DK5, Roeger3, Roeger4, Roeger5, Wood} \ldots . One popular approach to the elementary stability is  the investigation of Jacobian matrices of the discrete models at the equilibria, namely, determination of conditions ensuring that all eigenvalues of Jacobian matrices have moduli less or equal 1. This is the necessary and sufficient condition for hyperbolic equilibrium points to be locally stable \cite{Linda, Keshet}. The mentioned above approach has the following weaknesses and limitations:\par
\begin{enumerate}
\item It is applicable when all the equilibrium points are hyperbolic. To our best knowledge, at present no results on NSFD schemes preserving the stability of non-hyperbolic equilibrium points are available.
\item The consideration of Jacobian only guarantees the local stability meanwhile  many models have the global stability.
\end{enumerate}
Return to the model \eqref{eq:1}.  A difficulty is that when  $m_1 = r(0)$ or $m_2 = s(0)$ the equilibrium point $P_0^* = (0, 0)$ becomes non-hyperbolic. Therefore, it is impossible to study its stability via the eigenvalues of Jacobian $J(P_0^*)$. Consequently, we cannot directly apply the results concerning the elementary stability \cite{DK1, DK2, K1} of the system. For the continuous model, the Lyapunov stability theorem should be used for investigating the global stability \cite[Theorem 1]{Lindano}.\par

The models considered in  \cite{Bairagi, Darti, DK3, DK5} have equilibrium points, which are all hyperbolic. Meanwhile the system \eqref{eq:1} has a non-hyperbolic equilibrium point in a particular case of the parameters. Moreover, it is globally asymptotical stable.  Therefore, for the corresponding discrete system we shall use the Lyapunov stability theorem for proving the stability of this equilibrium point. This is the important contribution of our paper. Besides, we propose a more general model with many iterative parameters in the discretization of the right-hand sides. The combination of appropriate selection of denominator with  these parameters will give sufficient conditions for dynamical consistency of the discrete model with the continuous system.
 
\section{Construction of NSFD scheme}
In this section we construct NSFD scheme for the system \eqref{eq:1}  preserving all dynamic properties of the original continuous model for any discretization parameter or step size $h > 0$. Recall that, according to Mickens, a finite difference scheme is called \emph{nonstandard} if at least one of the following conditions is satisfied 
 \cite{Mickens1, Mickens2, Mickens3, Mickens4}:
\begin{itemize}
\item A nonlocal approximation is used.
\item  The discretization of the derivative is not traditional and uses a function $0 < \varphi(h) = h + \mathcal{O}(h^2)$.
\end{itemize}
We propose a general NSFD scheme for the model \eqref{eq:1} in the form
\begin{equation}\label{eq:3}
\begin{split}
\dfrac{x_{k + 1} - x_k}{\varphi(h)} &= \alpha_1 x_k r(x_k) + \alpha _2x_{k + 1}r(x_k) - \alpha_3 x_k y_k \phi(x_k) - \alpha_4x_{k + 1}y_k\phi(x_k) - \alpha_5m_1x_k - \alpha_6 m_1 x_{k + 1},\\
\dfrac{y_{k + 1} - y_k}{\varphi(h)} &= \beta_1 y_k s(y_k) + \beta_2y_{k + 1}s(y_k) + c\beta_3x_ky_k\phi(x_k) + c\beta_4x_ky_{k + 1}\phi(x_k) - \beta_5 m_2y_k - \beta_6 m_2 y_{k + 1},\\
\alpha_j + \alpha_{j + 1} &= \beta_j + \beta_{j + 1} = 1, \quad j = 1, 3, 5; \quad \varphi = h + \mathcal{O}(h^2), \quad h \to 0.
\end{split}
\end{equation}
Concerning the set of equilibrium points and positivity of the scheme \eqref{eq:3} there hold the following assertions.
%The proof of the proposition is like to that of Proposition 1 in \cite{Lindano}.
\begin{proposition}\label{Proposition2}
The region $\Omega = \big\{(x, y) \in \mathbb{R}^2\big| x \geq 0, \quad y \geq 0\big\}$ is a positive invariant set for the scheme \eqref{eq:3} if
\begin{equation}\label{eq:4}
\begin{split}
&\alpha_1 \geq 0, \quad \alpha_2 \leq 0, \quad \alpha_3 \leq 0, \quad \alpha_4 \geq 0, \quad \alpha_5 \leq 0, \quad \alpha_6 \geq 0,\\
&\beta_1 \geq 0, \quad \beta_2 \leq 0, \quad \beta_3 \geq 0, \quad \beta_4 \leq 0, \quad \beta_5 \leq 0, \quad \beta_6 \geq 0.
\end{split}
\end{equation}
\end{proposition}
\begin{proof}
It is easy to reduce the scheme \eqref{eq:3} to the explicit form
\begin{equation}\label{eq:3a}
\begin{split}
x_{k + 1} &= F(x_k, y_k) := \dfrac{x_k + \varphi \alpha_1 x_k r(x_k) - \varphi \alpha_3 x_k y_k \phi(x_k) - \varphi \alpha_5 m_1 x_k}{1 - \varphi \alpha_2 r(x_k) + \varphi \alpha_4 y_k \phi(x_k) + \varphi \alpha_6 m_1},\\
y_{k + 1}& = G(x_k, y_k):= \dfrac{y_k + \varphi \beta_1 y_k s(y_k) + \varphi \beta_3 c x_k y_k \phi(x_k) - \varphi \beta_5 m_2 y_k}{1 - \varphi \beta_2 s(y_k) - \varphi \beta_4 c x_k \phi(x_k) + \varphi \beta_6 m_2}.
\end{split}
\end{equation}
Since the parameters $\alpha_j$ and $\beta_j$ ($j = \overline{1, 6}$) satisfy \eqref{eq:4} from the above formulas we obtain the fact to be proved.
\end{proof}
\begin{proposition}\label{Propositon1}
 The scheme \eqref{eq:3} preserves the set of equilibrium points of system \eqref{eq:3}.
\end{proposition}
\begin{proof}
The equilibrium $(x^*, y^*)$ of the scheme \eqref{eq:3} is the solution of the equation
\begin{equation*}
x^* = F(x^*, y^*), \qquad y^* = G(x^*, y^*), 
\end{equation*}
where $F(x, y)$ and $G(x, y)$ are defined by \eqref{eq:3a}. Since the parameters $\alpha_i, \beta_i$ satisfy \eqref{eq:3}, if the scheme \eqref{eq:3} is defined then it is equivalent to
\begin{equation*}
x^*\big[r(x^*) - y^*\phi(x^*) - m_1 \big] = 0,\qquad y^*\big[s(y^*) + cx^*\phi(x^*) - m_2\big] = 0.
\end{equation*}
It is the system for determining equilibrium points of the model \eqref{eq:1}. Thus, the proposition is proved.
\end{proof}
\section{Stability analysis}
In this section we give sufficient conditions for the scheme \eqref{eq:3} to preserve the stability properties of equilibrium points of the model \eqref{eq:1}. For the purpose of easy tracking we recall the following results of the necessary and sufficient conditions for an equilibrium point to be locally asymptotically stable \cite[Theorem 2.10]{Linda}, \cite{Keshet}.
\begin{lemma}\label{Lemma1}
Assume the functions $f(x, y)$ and $g(x, y)$ have continuous first-order partial
derivatives in $x$ and $y$ on some open set in $\mathbb{R}^2$ that contains the point $(x^*, y^*)$. Then the equilibrium point $(x^*, y^*)$ of the nonlinear system 
\begin{equation*}
x_{k + 1} = f(x_k, y_k), \qquad y_{k + 1} = g(x_k, y_k),
\end{equation*}
is locally asymptotically stable if the eigenvalues of the Jacobian matrix J evaluated
at the equilibrium satisfy $|\lambda_i| < 1$ iff
\begin{equation*}
(i) \det(J) < 1, \qquad (ii) 1 - Tr(J) + det(J) > 0, \qquad (iii) 1 + Tr(J) + \det(J) > 0. 
\end{equation*}
The equilibrium is unstable if some $|\lambda_i| > 1$, that is, if any one of three inequalities is satisfied,
\begin{equation*}
(i) \det(J) > 1, \qquad (ii) 1 - Tr(J) + det(J) < 0, \qquad (iii) 1 + Tr(J) + \det(J) < 0. 
\end{equation*}
\end{lemma}
\subsection{The extinction equilibrium point $P_0^* = (0, 0)$}
\begin{proposition}\label{Proposition1}
For the case $m_1 > r(0)$ and $m_2 > s(0)$ consider the difference scheme \eqref{eq:3} under the assumptions of Proposition \ref{Proposition2}.
Then the extinction equilibrium point $P_0^* = (x_0^*, y_0^*) = (0, 0)$ is locally asymptotically stable if $m_1 > r(0)$ and $m_2 > s(0)$, and unstable otherwise.
\end{proposition}

\begin{proof}
Computing the Jacobian matrix of system \eqref{eq:3}, evaluated at the extinction equilibrium point $P_0^* = (x_0^*, y_0^*) = (0, 0)$ , one obtains
\begin{equation*}
J(P_0^*) = 
\begin{pmatrix}
\dfrac{1 + \varphi \alpha_1 r(0) - \varphi \alpha_5 m_1}{1 - \varphi \alpha_2 r(0) +\varphi m_1 \alpha_6}& 0\\
0 & \dfrac{1 + \varphi \beta_1 s(0) - \varphi \beta_5 m_2}{1 - \varphi \beta_2 s(0) + \varphi \beta_6 m_2}
\end{pmatrix}.
\end{equation*}
In this case, the eigenvalues are
\begin{equation*}
\lambda_1 = \dfrac{1 + \varphi \alpha_1 r(0) - \varphi \alpha_5 m_1}{1 - \varphi \alpha_2 r(0) +\varphi m_1 \alpha_6}, \qquad \lambda_2 = \dfrac{1 + \varphi \beta_1 s(0) - \varphi \beta_5 m_2}{1 - \varphi \beta_2 s(0) + \varphi \beta_6 m_2}.
\end{equation*}
Since  $\alpha_j$ and $\beta_j$ satisfy Proposition \ref{Proposition2} we have $\lambda_1 > 0, \, \, \lambda_2 > 0$. On the other hand
\begin{equation*}
|\lambda_1| - 1 = \varphi \dfrac{r(0) - m_1}{1 - \varphi \alpha_2 r(0) +\varphi m_1 \alpha_6}, \qquad |\lambda_2| - 1 = \varphi \dfrac{s(0) - m_2}{1 - \varphi \beta_2 s(0) + \varphi \beta_6 m_2},
\end{equation*}
Therefore, $ |\lambda_1|, |\lambda_2| < 1$ iff $m_1 > r(0)$ and $m_2 > s(0)$. By Lemma \ref{Lemma1} Proposition is proved.
\end{proof}
Next, consider scheme \eqref{eq:3} under the assumptions of  Proposition \ref{Proposition2} in the case $m_1 \geq r(0)$ and $m_2 \geq s(0)$. Notice that when $m_1 = r(0)$ or $m_2 = s(0)$ the equilibrium point $P_0^* = (0, 0)$ becomes non-hyperbolic one. Therefore, it is impossible to use  Lemma \ref{Lemma1} for proving the locally asymptotical stability. Moreover, the equilibrium point $P_0^*$, as shown in \cite[Theorem 1]{Lindano}, is globally asymptotically stable. As for the continuous system we shall also use Lyapunov's stability theorem to prove that the point $P_0^* = (0, 0)$ is a globally asymptotically stable equilibrium point of the scheme \eqref{eq:3}. For this purpose we consider a family of functions
\begin{equation}\label{eq:6}
V(x_k, y_k) := \alpha x_k y_k + \beta x_k^2 + \gamma x_k + \delta y_k, \quad (x_k, y_k) \in \mathbb{R}_+^2,
\end{equation}
where $\alpha, \beta, \gamma, \delta > 0$ are parameters, which are selected so that the function $V(x_k, y_k)$ satisfies  the conditions of the Lyapunov's stability theorem \cite[Theorem 4.20]{Elaydi}.\par
Obviously, the function $V(x_k, y_k)$ defined by \eqref {eq:6} is continuous on $\mathbb{R}_+^2$ and $V(x_k, y_k) \to \infty$ as $||(x_y, y_k)|| \to \infty$. Moreover, $V(P_0^*) = 0$ and $V(x_k, y_k) > 0$ for any $(x_k, y_k) \in \mathbb{R}_+^2$, $(x_k, y_k) \ne (0, 0)$. Therefore, in order to show that $V(x_k, y_k)$ satisfies \cite[Theorem 4.20]{Elaydi} it  suffices to determine conditions for 
%ta chỉ xác định các điều kiện để
\begin{equation*}
 \Delta V(x_k, y_k) = V(x_{k + 1}, y_{k + 1}) - V(x_k, y_k) < 0  \quad \forall(x_k, y_k) \in \mathbb{R}_+^2\backslash\{(0, 0)\}.
\end{equation*}
We have
\begin{equation}\label{eq:7}
\Delta V(x_k, y_k) = \alpha (x_{k + 1} y_{k + 1} - x_ky_k) + \beta (x_{k + 1}^2 - x_k^2) + \gamma (x_{k + 1} - x_k) + \delta (y_{k + 1} - y_k).
\end{equation}
Now we rewrite \eqref{eq:3a} in the form
\begin{equation}\label{eq:3b}
\begin{split}
x_{k + 1} &=  x_k + \varphi \dfrac{x_k[r(x_k) - m_1] - x_k y_k\phi(x_k)}{1 - \varphi \alpha_2 r(x_k) + \varphi \alpha_4 y_k \phi(x_k) + \varphi \alpha_6 m_1},\\
y_{k + 1} &= y_k + \varphi\dfrac{y_k[s(y_k) - m_2] + cx_ky_k\phi(x_k)}{{1 - \varphi \beta_2 s(y_k) - \varphi \beta_4 c x_k \phi(x_k) + \varphi \beta_6 m_2}}.
\end{split}
\end{equation}
According to the estimates in Theorem 1 in \cite{Lindano} we have: if $m_1 \geq r(0)$ and $m_2 \geq s(0)$, then $r(x) - m_1 \leq 0$ for all $x \geq 0$ and $s(y) - m_2 \leq 0$ for all $y \geq 0$. Therefore, from \eqref{eq:3b} we see that if $m_1 \geq r(0)$, then $x_{k + 1} \leq x_k$ for all $k \geq 0$. Consequently, from \eqref{eq:7} it implies that if $y_{k + 1} - y_k < 0$ then $\Delta V(x_k, y_k) < 0$. Hence, it is sufficient to consider only the case $y_{k + 1} \geq y_k$.\par
Using the mean value theorem for two variables function $u(x_k, y_k) = x_k y_k$ we have
\begin{equation}\label{eq:8}
x_{k + 1}y_{k + 1} - x_k y_k = \xi_{y_k}(x_{k + 1} - x_k) + \xi_{x_k}(y_{k + 1} - y_k) \leq y_k (x_{k + 1} - x_k) + x_k(y_{k + 1} - y_k),
\end{equation}
where $\xi_{x_k} \in (x_{k + 1}, x_k)$ and $\xi_{y_k} \in (y_k, y_{k + 1})$. On the other hand, since $x_{k + 1} \leq x_k$ we have
\begin{equation}\label{eq:9}
x_{k + 1}^2 - x_{k}^2 \leq x_kx_{k + 1} - x_k^2 = x_k(x_{k + 1} - x_k).
\end{equation}
From \eqref{eq:8}, \eqref{eq:9}, \eqref{eq:7} and \eqref{eq:3b} we obtain the estimate
\begin{equation}\label{eq:10}
\begin{split}
\Delta V(x_k, y_k) &\leq (\alpha y_k +  \beta x_k + \gamma)(x_{k + 1} - x_k) + (\alpha x_k + \delta)(y_{k + 1} - y_k)\\
&= (\alpha y_k + \beta x_k + \gamma) \varphi \dfrac{x_k[r(x_k) - m_1] - x_k y_k\phi(x_k)}{1 - \varphi \alpha_2 r(x_k) + \varphi \alpha_4 y_k \phi(x_k) + \varphi \alpha_6 m_1}\\
 &+  (\alpha x_k + \delta) \varphi \dfrac{y_k[s(y_k) - m_2] + cx_ky_k\phi(x_k)}{{1 - \varphi \beta_2 s(y_k) - \varphi \beta_4 c x_k \phi(x_k) + \varphi \beta_6 m_2}}\\
&\leq \varphi x_ky_k\phi(x_k) Q(x_k, y_k),
%F(x_k, y_k, \varphi, \alpha, \beta, \gamma, \delta ),
% \big[big]
\end{split}
\end{equation}
where $Q(x_k, y_k)$ is defined by
\begin{equation*}
Q(x_k, y_k) := \dfrac{c(\alpha x_k +\delta)}{1 - \varphi \beta_2 s(y_k) - \varphi \beta_4 c x_k \phi(x_k) + \varphi \beta_6 m_2} - \dfrac{\alpha y_k +  \beta x_k + \gamma}{1 - \varphi \alpha_2 r(x_k) + \varphi \alpha_4 y_k \phi(x_k) + \varphi \alpha_6 m_1}.
\end{equation*}
From here it follows that if $Q(x_k, y_k) < 0$ for all $(x_k, y_k) \in \mathbb{R}_+^2\backslash\{(0, 0)\}$ then $\Delta V(x_k, y_k) < 0$ for all $(x_k, y_k) \in \mathbb{R}_+^2\backslash\{(0, 0)\}$. On the other hand we have
\begin{equation*}
\begin{split}
Q(x_k, y_k)[1 - \varphi \alpha_2 r(x_k) + \varphi \alpha_4 y_k \phi(x_k) + \varphi \alpha_6 m_1][{1 - \varphi \beta_2 s(y_k) - \varphi \beta_4 c x_k \phi(x_k) + \varphi \beta_6 m_2}]\\
= \tau_1 x_k + \tau_2 \varphi x_k + \tau_3 \varphi c \alpha x_k y_k \phi(x_k) + \tau_4 + \varphi \tau_5 + \tau_6 \varphi y_k - Q_1(x_k, y_k),
\end{split}
\end{equation*}
where
\begin{equation}\label{eq:tau_i}
\begin{split}
&\tau_1 := c \alpha - \beta, \quad \tau_2 := -\alpha_2 c \alpha r(x_k) + \alpha_6 c \alpha m_1 - \beta \beta_6 m_2, \qquad \tau_3 := \alpha_4 + \beta_4,\\
&\tau_4 := c\delta - \gamma, \quad \tau_5 := -\alpha_2 c \delta r(x_k) + c \delta \alpha_6 m_1 - \beta_6 \gamma m_2, \quad \tau_6 := \alpha_4 c \delta \phi(x_k) - \alpha \beta_6 m_2,\\
&Q_1(x_k, y_k) := \beta x_k [-\varphi \beta_2 s(y_k) - \varphi \beta_4 c x_k \phi(x_k)] + \alpha y_k[1 - \varphi \beta_2 s(y_k)] + \gamma[-\varphi \beta_2 s(y_k) - \varphi \beta_4 c x_k \phi(x_k)].
\end{split}
\end{equation}
Notice that $Q_1(x_k, y_k) > 0$ for all  $(x_k, y_k) \in \mathbb{R}_+^2\backslash\{(0, 0)\}$. Therefore, $Q(x_k, y_k) < 0$ for all  $(x_k, y_k) \in \mathbb{R}_+^2\backslash\{(0, 0)\}$ if $\tau_i < 0$ \, ($i = \overline{1, 6}$). Taking into account \eqref{eq:2} we have $0 < r(x_k) \leq r(0)$ and $0 <\phi(x_k) \leq \phi(0)$ for all $x_k \geq 0$. Therefore from \eqref{eq:tau_i} it follows
\begin{enumerate}
\item $\tau_1 < 0$ if $c < \dfrac{\beta}{\alpha}$,
\item $\tau_2 = -\alpha_2 c \alpha r(x_k) + \alpha_6 c \alpha m_1 - \beta \beta_6 m_2$ $\leq$  $-\alpha_2 c \alpha r(0) + \alpha_6 c \alpha m_1 - \beta \beta_6 m_2 < 0$\\
 if \quad $\dfrac{-\alpha_2 c r(0) + \alpha_6 c m_1}{\beta_6 m_2} < \dfrac{\beta}{\alpha}$,
\item $\tau_3  < 0$ if $\alpha_4 + \beta_4 < 0$,
\item $\tau_4  < 0$ if $c < \dfrac{\gamma}{\delta}$,
\item $\tau_5 = -\alpha_2 c \delta r(x_k) + c \delta \alpha_6 m_1 - \beta_6 \gamma m_2 \leq -\alpha_2 c \delta r(0) + c \delta \alpha_6 m_1 - \beta_6 \gamma m_2 < 0$ \\
if \quad $\dfrac{-\alpha_2 c r(0) + c \alpha_6 m_1}{\beta_6 m_2} < \dfrac{\gamma}{\delta}$,
\item $\tau_6 = \alpha_4 c \delta \phi(x_k) - \alpha \beta_6 m_2 \leq \alpha_4 c \delta \phi(0) - \alpha \beta_6 m_2 < 0$\\
 if \quad $\dfrac{c \alpha_4 \phi(0)}{\beta_6m_2} < \dfrac{\alpha}{\delta}$.
\end{enumerate}
In summary the function $V(x_k, y_k)$ defined by \eqref{eq:6} satisfies $ \Delta V(x_k, y_k) < 0$  for all $(x_k, y_k) \in \mathbb{R}_+^2\backslash\{(0, 0)\}$ if
\begin{equation}\label{eq:11}
\begin{split}
&\max\Big\{c,\dfrac{-\alpha_2 c r(0) + \alpha_6 c m_1}{\beta_6 m_2}\Big\} < \dfrac{\beta}{\alpha}, \quad \max\Big\{c,\dfrac{-\alpha_2 c r(0) + c \alpha_6 m_1}{\beta_6 m_2}\Big\} < \dfrac{\gamma}{\delta},\\
& \dfrac{c \alpha_4 \phi(0)}{\beta_6m_2} < \dfrac{\alpha}{\delta}, \qquad \alpha_4 + \beta_4 < 0.
\end{split}
\end{equation}
Once, the scheme \eqref{eq:3} is fixed the selection of the parameters $\alpha, \beta, \gamma, \delta > 0$ satisfying the above relations is completely possible.
Thus, we obtain the following theorem of the global stability of the equilibrium point $P_0^* = (0, 0)$.
\begin{theorem}\label{Proposition3}
For the case $m_1 \geq r(0)$ and $m_2 \geq s(0)$ consider the difference scheme \eqref{eq:3} under the assumptions of Proposition \ref{Proposition2}. If additionally assume that 
\begin{equation}\label{eq:p1}
\alpha_4 + \beta_4 < 0,
\end{equation}
 then the extinction equilibrium point $P_0^* = (0, 0)$ is globally asymptotically stable.
\end{theorem}
\subsection{The  equilibrium point $P_1^* = (K, 0)$ (the equilibrium point of extinction of the predator species)}
\begin{proposition}\label{Proposition4}
For the case $m_1 < r(0)$ and $m_2 > s(0) + cK\phi(K)$ consider the difference scheme \eqref{eq:3} under the assumptions of Proposition \ref{Proposition2}.
 If additionally assume that
\begin{equation}\label{eq:12}
\begin{split}
&T_1 := 2 \alpha_6 m_1 - 2 \alpha_2 r(K) + K r'(K) > 0,\\
&T_2 := s(0) - m_2 + cK\phi(K) - 2 \beta_2 s(0) - 2 \beta_4 c K \phi(K) + 2 \beta_6 m_2 > 0,
\end{split}
\end{equation}
then the equilibrium point of the form $P_1^* = (K, 0)$ is locally asymptotically stable if  $m_1 < r(0)$ and $m_2 > s(0) + cK\phi(K)$, and unstable otherwise.
\end{proposition}
\begin{proof}
Recall that $m_1 < r(0)$ is necessary and sufficient condition for the existence of the equilibrium point $P_1^*$ (see Theorem 1).
Computing the Jacobian matrix of system \eqref{eq:3}, evaluated at the extinction equilibrium point $P_1^* = (K, 0)$, one obtains
\begin{equation*}
J(P_1^*) = 
\begin{pmatrix}
1 + \dfrac{\varphi  K r'(K)}{1 - \varphi \alpha_2 r(K) + \varphi \alpha_6 m_1}& \dfrac{-\varphi K \phi(K)}{1 - \varphi \alpha_2 r(K) + \varphi \alpha_6 m_1}\\
0 & 1 + \varphi \dfrac{s(0) - m_2 + cK\phi(K)}{1 - \varphi \beta_2 s(0) - \varphi \beta_4 c K \phi(K) + \varphi \beta_6 m_2}
\end{pmatrix}.
\end{equation*}
In this case, the eigenvalues are
\begin{equation*}
\lambda_1 = 1 + \dfrac{\varphi  K r'(K)}{1 - \varphi \alpha_2 r(K) + \varphi \alpha_6 m_1}, \qquad \lambda_2 =  1 + \varphi \dfrac{s(0) - m_2 + cK\phi(K)}{1 - \varphi \beta_2 s(0) - \varphi \beta_4 c K \phi(K) + \varphi \beta_6 m_2}.
\end{equation*}
Therefore, $\lambda_2 < 1$ iff  $m_2 > s(0) + c K \phi(K)$. Besides, due to $r'(K) < 0$  there holds $\lambda_1 < 1$. On the other hand we have

\begin{equation*}
% \begin{split}
\lambda_1 + 1 = \dfrac{2 + \varphi T_1}{1 - \varphi \alpha_2 r(K) + \varphi \alpha_6 m_1},\quad \lambda_2 + 1 = \dfrac{2 + \varphi T_2}{1 - \varphi \beta_2 s(0) - \varphi \beta_4 c K \phi(K) + \varphi \beta_6 m_2},
% \end{split}
\end{equation*}
Therefore, if  \eqref{eq:12} is fulfilled then $\lambda_1 > - 1$, $\lambda_2 > - 1$. In result we have $|\lambda_1| < 1$ and $|\lambda_2| < 1$. Hence, by Lemma \ref{Lemma1} the proposition is proved.
\end{proof}

%%%%%%%%%%%%%%%%%%%%%%%%%%%%%%%%%%%%%%%%%%%%%%%%%%%%%%%%%%%%%%%%%%%%%%%%%%%%%%%%%%%%%%%%%%%%%%%%%%%%%%%%%%%%%
\subsection{The  equilibrium point $P_2^* = (0, M)$ (the equilibrium point of extinction of the prey species)}
\begin{proposition}\label{Proposition5}
For the case $m_1 > r(0) - M\phi(0)$ and $m_2 < s(0)$ consider the difference scheme \eqref{eq:3} under the assumptions of Proposition \ref{Proposition2}.
 If additionally assume that
\begin{equation}\label{eq:13}
\begin{split}
& T_3 := r(0) - M \phi (0) - m_1 - 2 \alpha_2 r(0) + 2 \alpha_4 M \phi(0) + 2 \alpha_6 m_1 > 0,\\
&T_4 := M s'(M) - 2 \beta_2 s(M) + 2 \beta_6 m_2 >0,
\end{split}
\end{equation}
then the equilibrium point of the form $P_2^* = (0, M)$ is locally asymptotically stable if  $m_1 > r(0) - M\phi(0)$ and $m_2 < s(0)$, and unstable otherwise.
\end{proposition}
\begin{proof}
Notice that according to Theorem 1, $m_2 < s(0)$ is necessary and sufficient condition for the existence of the equilibrium point $P_2^*$. Computing the Jacobian matrix of system \eqref{eq:3}, evaluated at the extinction equilibrium point $P_2^* = (0, M)$, one obtains
\begin{equation*}
J(P_2^*) = 
\begin{pmatrix}
1 + \varphi \dfrac{r(0) - M \phi (0) - m_1}{1 - \varphi \alpha_2 r(0) + \varphi \alpha_4 M \phi(0) + \varphi \alpha_6 m_1}& 0\\
\dfrac{\varphi c M \phi(0)}{1 - \varphi \beta_2 s(M) + \varphi \beta_6 m_2} & 1 +  \dfrac{\varphi M s'(M)}{1 - \varphi \beta_2 s(M) + \varphi \beta_6 m_2}
\end{pmatrix}.
\end{equation*}
In this case, the eigenvalues are
\begin{equation*}
\lambda_1 = 1 + \varphi \dfrac{r(0) - M \phi (0) - m_1}{1 - \varphi \alpha_2 r(0) + \varphi \alpha_4 M \phi(0) + \varphi \alpha_6 m_1}, \quad \lambda_2 =  1 +  \dfrac{\varphi M s'(M)}{1 - \varphi \beta_2 s(M) + \varphi \beta_6 m_2}.
\end{equation*}
Therefore, $\lambda_1 < 1$ iff $m_1 > r(0) - M\phi(0)$. Besides, due to $s'(M) < 0$ there holds $\lambda_2 < 1$. On the other hand we have
\begin{equation*}
\begin{split}
\lambda_1 + 1 = \dfrac{2 +\varphi T_3}{1 - \varphi \alpha_2 r(0) + \varphi \alpha_4 M \phi(0) + \varphi \alpha_6 m_1}, \quad \lambda_2 + 1 = \dfrac{2 + \varphi T_4}{1 - \varphi \beta_2 s(M) + \varphi \beta_6 m_2}.
\end{split}
\end{equation*}
Hence, if \eqref{eq:13} is fulfilled then  $\lambda_2 > -1$ and $\lambda_1 > -1$. In result we have $|\lambda_1| < 1$ and $|\lambda_2| < 1$. By Lemma \ref{Lemma1} the  proposition is proved.
\end{proof}
%%%%%%%%%%%%%%%%%%%%%%%%%%%%%%%%%%%%%%%%%%%%%%%%%%%%%%%%%%%%%%%%%%%%%%%%%%%%%%%%%%%%%%%%%%%%%%%%%%%%%%%%%%%%%%%%%%%%%%%%%%%%%%%%%%%%%%%%%%%%%%%%%%%%%%%%%%%%%%%%%%%%%%%%%%%
\subsection{The  equilibrium point $P_3^* = (x^*, y^*)$ (the ecological stability equilibrium)}
\begin{proposition}\label{Proposition4}
Consider the difference scheme \eqref{eq:3} under the assumptions of Proposition \ref{Proposition2}.
 If additionally assume that
\begin{equation}\label{eq:14}
\begin{split}
T_5 := &-x^*[r'(x^*) - y^*\phi'(x^*)][- \beta_2 s(y^*) - \beta_4 c x^* \phi(x^*) + \beta_6 m_2]\\
& - y^*s'(y^*)[- \alpha_2 r(x^*) + \alpha_4 y^* \phi(x^*) + \alpha_6 m_1]\\
& -  x^*y^*s'(y^*)[r'(x^*) - y^*\phi'(x^*)] - cx^*y^*\phi(x^*)[\phi(x^*) + x^*\phi'(x^*)] > 0,\\
T_6:=& -\alpha_2 r(x^*) + \alpha_4 y^* \phi(x^*) + \alpha_6 m_1  +  x^*[r'(x^*) - y^*\phi'(x^*)] > 0,\\
 T_7 :=& -\beta_2 s(y^*) - \beta_4 c x^* \phi(x^*) + \beta_6 m_2 + y^*s'(y^*) > 0, 
\end{split}
\end{equation}
then, if the equilibrium point of the form $P_3^* = (x^*, y^*)$ belong to $\Omega$, then it is locally asymptotically stable.
\end{proposition}
\begin{proof}
Notice that the conditions for the existence of the equilibrium point $P_3^*$ are given in Theorem 1.
Computing the Jacobian matrix of system \eqref{eq:3}, evaluated at the extinction equilibrium point $P_3^* = (x^*, y^*)$, one obtains
\begin{equation*}
J(P_3^*) = 
\begin{pmatrix}
1 +  \dfrac{\varphi x^*[r'(x^*) - y^*\phi'(x^*)]}{1 - \varphi \alpha_2 r(x^*) + \varphi \alpha_4 y^* \phi(x^*) + \varphi \alpha_6 m_1}& \dfrac{-\varphi x^* \phi(x^*)}{1 - \varphi \alpha_2 r(x^*) + \varphi \alpha_4 y^* \phi(x^*) + \varphi \alpha_6 m_1}\\
\\
\varphi \dfrac{cy^*[\phi(x^*) + x^*\phi'(x^*)]}{1 - \varphi \beta_2 s(y^*) - \varphi \beta_4 c x^* \phi(x^*) + \varphi \beta_6 m_2} & 1 +  \dfrac{\varphi y^*s'(y^*)}{1 - \varphi \beta_2 s(y^*) - \varphi \beta_4 c x^* \phi(x^*) + \varphi \beta_6 m_2}
\end{pmatrix}.
\end{equation*}
In this case, we have
\begin{equation*}
\begin{split}
\det(J(P_3^*)) &= 1 + \dfrac{\varphi x^*[r'(x^*) - y^*\phi'(x^*)]}{u(x^*, y^*)} + \dfrac{\varphi y^*s'(y^*)}{v(x^*, y^*)}\\
 &+  \varphi^2\dfrac{x^*y^*s'(y^*)[r'(x^*) - y^*\phi'(x^*)] + cx^*y^*\phi(x^*)[\phi(x^*) + x^*\phi'(x^*)]}{u(x^*, y^*)v(x^*, y^*)},\\
Tr(J(P_3^*)) &= 2 +  \dfrac{\varphi x^*[r'(x^*) - y^*\phi'(x^*)]}{u(x^*, y^*)} + \dfrac{\varphi y^*s'(y^*)}{v(x^*, y^*)},
\end{split}
\end{equation*}
where
\begin{equation*}
\begin{split}
u(x^*, y^*) & := 1 - \varphi \alpha_2 r(x^*) + \varphi \alpha_4 y^* \phi(x^*) + \varphi \alpha_6 m_1 > 0,\\
v(x^*, y^*) & := 1 - \varphi \beta_2 s(y^*) - \varphi \beta_4 c x^* \phi(x^*) + \varphi \beta_6 m_2 > 0.
\end{split}
\end{equation*}
According to Theorem 2 in \cite{Lindano} we have $s'(y^*) < 0,\,r'(x^*) - y^*\phi'(x^*) < 0,$  consequently
\begin{equation*}
\begin{split}
1 - Tr(J(P_3^*)) + \det(J(P_3^*)) =
 \varphi^2\dfrac{x^*y^*s'(y^*)[r'(x^*) - y^*\phi'(x^*)] + cx^*y^*\phi(x^*)[\phi(x^*) + x^*\phi'(x^*)]}{u(x^*, y^*)v(x^*, y^*)} > 0.
\end{split}
\end{equation*}
On the other hand  $\det(J(P_3^*)) < 1$ iff
\begin{equation*}
-\dfrac{ x^*[r'(x^*) - y^*\phi'(x^*)]}{u(x^*, y^*)} - \dfrac{y^*s'(y^*)}{v(x^*, y^*)}
 >  \varphi\dfrac{x^*y^*s'(y^*)[r'(x^*) - y^*\phi'(x^*)] + cx^*y^*\phi(x^*)[\phi(x^*) + x^*\phi'(x^*)]}{u(x^*, y^*)v(x^*, y^*)}.
\end{equation*}
This equivalent to
\begin{equation*}
\dfrac{-x^*[r'(x^*) - y^*\phi'(x^*)] - y^*s'(y^*)}{u(x^*, y^*)v(x^*, y^*)} + \dfrac{\varphi  T_5}{u(x^*, y^*)v(x^*, y^*)} > 0.
\end{equation*}
The first term of the left-hand side of the above inequality is positive, therefore, if $T_5 > 0$ then the above  inequality is valid. It implies $\det(J(P_3^*)) < 1$. \par
Finally, we see that $1 + Tr(J(P_3^*)) + \det (J(P_3^*)) > 0$ iff
\begin{equation*}
\begin{split}
2\Big[1 &+ \dfrac{\varphi x^*[r'(x^*) - y^*\phi'(x^*)]}{u(x^*, y^*)}\Big] + 2\Big[1 + \dfrac{\varphi y^*s'(y^*)}{v(x^*, y^*)}\Big]\\
 &+  \varphi^2\dfrac{x^*y^*s'(y^*)[r'(x^*) - y^*\phi'(x^*)] + cx^*y^*\phi(x^*)[\phi(x^*) + x^*\phi'(x^*)]}{u(x^*, y^*)v(x^*, y^*)} > 0
\end{split}
\end{equation*}
The third term of the above sum is always positive, so if
\begin{equation}\label{eq:15}
1 + \dfrac{\varphi x^*[r'(x^*) - y^*\phi'(x^*)]}{u(x^*, y^*)} > 0, \qquad 
1 + \dfrac{\varphi y^*s'(y^*)}{v(x^*, y^*)} > 0,
\end{equation}
then $1 + Tr(J(P_3^*)) + \det (J(P_3^*)) > 0$. It is easy to verify that if $T_6, T_7 > 0$ then \eqref{eq:15} holds. Thus, we have proved that if 
 \eqref{eq:14} satisfies then the three conditions of Lemma \ref{Lemma1} are satisfied. Therefore, the point $P_3^*$  is locally asymptotically stable. The proposition is proved.
%\begin{equation*}
%\det (J(P_3^*)) < 1, \quad 1 - trace(J(P_3^*)) + \det (J(P_3^*)) > 0, \quad 1 + trace(J(P_3^*)) + \det (J(P_3^*)) > 0,
%\end{equation*}
\end{proof}

\begin{remark}\label{Remark1}
The system of conditions \eqref{eq:p1}-\eqref{eq:14} for $\alpha_j$, $\beta_j$ \, $(j = \overline{1, 6})$ has many solutions. Simply, we can choose $\alpha_4, \alpha_6, \beta_6 > 0$ sufficiently large, $\alpha_2, \beta_2, \beta_4 < 0$  sufficiently small and $\alpha_4 + \beta_4 < 0$.
\end{remark}
Now we summarize the obtained results above in the following theorem on the NSFD schemes preserving the dynamical properties of the model \eqref{eq:1}.
\begin{theorem}\label{Maintheorem}
The NSFD scheme \eqref{eq:3} is dynamically consistent with \eqref{eq:1} if the parameters $\alpha_j$, $\beta_j$ \, $(j = \overline{1, 6})$ satisfy the conditions listed in Table \ref{tabl1}, where the columns list sufficient conditions for the scheme \eqref{eq:3} preserve corresponding properties of the model \eqref{eq:1} for different cases of the parameters. The symbol  $''*''$ means that the set of equilibrium points of \eqref{eq:1} is always preserved by the scheme  \eqref{eq:3}.
\begin{table}
\setlength{\tabcolsep}{0.12cm}
\caption{The sufficient conditions for dynamical consistency}\label{tabl1}
\medbreak
\begin{tabular}{ l  c  c  c }
\hline
\\
$(m_1, m_2)$& Set of equilibria & Positivity & Stability\\
\\
\hline 
\\
$m_1 \geq r(0)$ and $m_2 \geq s(0)$&*&\eqref{eq:4}&\eqref{eq:p1}\\ 
\\
\hline
\\
$m_1 < r(0)$ and $m_2 > s(0) + cK\phi(K)$&*&\eqref{eq:4}&\eqref{eq:12}\\ 
\\
\hline
\\
$m_1 > r(0) - M\phi(0)$ and $m_2 < s(0)$&*&\eqref{eq:4}&\eqref{eq:13}\\ 
\\
\hline
\\
$m_1 < r(0) - M\phi(0)$ and $m_2 < s(0)$&*&\eqref{eq:4}&\eqref{eq:14}\\ 
\\
\hline
\\
$m_1 < r(0)$ and $s(0) < m_2 < s(0) + cK\phi(K)$&*&\eqref{eq:4}&\eqref{eq:14}\\
\\
\hline
\end{tabular} 
\end{table}
\end{theorem}
\begin{remark}
There are infinitely many ways for selecting the parameters $\alpha_j$, $\beta_j$ \, $(j = \overline{1, 6})$ satisfying the conditions listed in Table \ref{tabl1} ( see  Remark \ref{Remark1}). This shows the existence of NSFD schemes dynamically consistent with the system   \eqref{eq:1}.
\end{remark}
\section{Numerical Simulations}
Numerical examples presented in this section show that the obtained theoretical results of the NSFD preserving the properties of the general predator-prey system are valid. \par
Let us consider the example of the predator-prey model presented in \cite{Lindano}. In this example
%Chúng tôi xem xét lại ví dụ được trình bày trong Section 4 in \cite{Lindano}. Trong ví dụ này 
\begin{equation*}
xr(x) = \dfrac{15x}{x + 10}, \qquad ys(y) = \dfrac{5y}{y + 10}, \qquad x\phi(x) = \dfrac{x}{x + 30}, \qquad c = 0.003,
\end{equation*}
for $6$ cases of the parameters $(m_1, m_2)$ in Corollary 1 in \cite{Lindano}. Namely,
\begin{equation*}
\begin{split}
&(i)\;  m_1 = 1.53, \quad m_2 = 0.622. \quad (ii)\;  m_1 = 1.53, \quad m_2 = 0.4789. \\
& (iii)\;  m_1 = 1.4925, \quad m_2 = 0.4789. \quad (iv)\;  m_1 = 1.38, \quad m_2 = 0.4789. \\
& (v)\; m_1 = 0.3, \quad m_2 = 0.501. \quad (vi)\;  m_1 = 1.38, \quad m_2 = 0.622.
\end{split}
\end{equation*}
Many numerical simulations, such as  \cite{AL1, AL2, DQA, DK1, DK2,DK3,DK4, Roeger4, Roeger5, Wood} \ldots agree that standard difference schemes do not preserve dynamical properties of continuous models for large step sizes, i.e., are not dynamically consistent with continuous systems. This  confirms the advantages of NSFD schemes. For supporting this confirmation in this example we use the explicit Euler scheme and four stage Runge-Kutta (RK4) scheme compared with the constructed NSFD schemes for the system \eqref{eq:1}. The numerical solutions obtained for these schemes in Case (i) of $m_1, \ m_2$ are depicted in Figures 1-4. From the figures we see that the property of positivity and stability of the 
system are destroyed. The numerical experiments in other cases of the parameters are analogous.\\
The numerical solutions obtained by these NSFD schemes are depicted in Figures 5-10, respectively. Comparing these results with the numerical simulations in Section 4 in \cite{Lindano} we see that all properties of the continuous model are preserved.\\
% \begin{figure}[!ht]
% \centering
% \includegraphics[height=8cm,width=15cm]{fig5.eps}
% \caption{The solution $(x_k, y_k)$ obtained by the scheme \eqref{eq:3}-(i) for $h = 10$, $t \in [0, 1000]$.}\label{fig:5}
% \end{figure}
% 
% 
% 
% \begin{figure}[!ht]
% \centering
% \includegraphics[height=8cm,width=15cm]{fig6.eps}
% \caption{The solution $(x_k, y_k)$ obtained by the scheme  \eqref{eq:3}-(ii) for $h = 5$, $t \in [0, 2000]$.}\label{fig:6}
% \end{figure}
% 
% 
% \begin{figure}[!ht]
% \centering
% \includegraphics[height=8cm,width=15cm]{fig7.eps}
% \caption{The solution $(x_k, y_k)$ obtained by the scheme  \eqref{eq:3}-(iii) for $h = 5$, $t \in [0, 2000]$.}\label{fig:7}
% \end{figure}
% 
% \begin{figure}[!ht]
% \centering
% \includegraphics[height=8cm,width=15cm]{fig8.eps}
% \caption{The solution $(x_k, y_k)$ obtained by the scheme  \eqref{eq:3}-(iv) for $h = 4$, $t \in [0, 2000]$.}\label{fig:8}
% \end{figure}
% 
% \begin{figure}[!ht]
% \centering
% \includegraphics[height=8cm,width=15cm]{fig9.eps}
% \caption{The solution $(x_k, y_k)$ obtained by the scheme \eqref{eq:3}-(v) for $h = 4$, $t \in [0, 2000]$.}\label{fig:9}
% \end{figure}
% \begin{figure}[!ht]
% \centering
% \includegraphics[height=8cm,width=15cm]{fig10.eps}
% \caption{The solution $(x_k, y_k)$ obtained by the scheme  \eqref{eq:3}-(vi) for $h = 5$, $t \in [0, 2000]$.}\label{fig:10}
% \end{figure}
% 
\section{Conclusion}
In this paper we have used  NSFD schemes for converting a general predator-prey model to a dynamically consistent discrete system. It should be emphasized that in the continuous as in the discrete system there is a non-hyperbolic equilibrium point, whose global asymptotical stability was proved by means of the Lyapunov stability theorem. The numerical simulations for the model considered in \cite{Lindano} for various collections of parameters confirm the validity of the obtained theoretical results. In the future we shall develop the techniques used in this paper for constructing and investigating NSFD schemes for other dynamical models including ones having non-hyperbolic equilibrium points.
\section*{Acknowledgments}
%The authors would like to thank the reviewers for their helpful comments and suggestions. \par
This work is supported by Vietnam National Foundation for Science and Technology Development (NAFOSTED) under the grant  number 102.01-2014.20.\\
 
\end{document}